\documentclass[ 12pt reqno,amsmath,amsthm,amssymb,amscd]{amsart}
\usepackage{mathrsfs,amssymb, amscd,amsmath,amsthm}
\usepackage[enableskew,vcentermath]{youngtab}
\usepackage{multicol}\multicolsep=0pt
\usepackage{pstricks,pst-node}
\usepackage[backref=page]{hyperref}
\renewcommand*{\backref}[1]{}
\renewcommand*{\backrefalt}[4]{{[\tiny%
  \ifcase #1\relax\or Page~#2.%
  \else Pages #2.\fi]}}
\usepackage[sort]{cite}

\def\crulefill{\leavevmode\leaders\hrule height 1pt\hfill\kern 0pt}
\long\def\QUERY#1{%
\leavevmode\newline%
\noindent$\star\star\star$\thinspace\textsf{Comment/Query}\crulefill\newline%
   \space #1\newline\hbox to 120mm{\crulefill}$\star\star\star$\newline
}

\numberwithin{equation}{section} \theoremstyle{definition}
\newtheorem{Defn}[equation]{Definition}%[section]
\newtheorem{Remark}[equation]{Remark}
\theoremstyle{plain}
\newtheorem{Prop}[equation]{Proposition}
\newtheorem{Theorem}[equation]{Theorem}

\newtheorem{Lemma}[equation]{Lemma}

\makeatletter
\def\enumerate{\begingroup\ifnum\@enumdepth>3\@toodeep\else
      \advance\@enumdepth\@ne
      \edef\@enumctr{enum\romannumeral\the\@enumdepth}%
      \topsep\z@\parskip\z@
      \list{\csname label\@enumctr\endcsname}
        {\@nmbrlisttrue\let\@listctr\@enumctr
         \parsep\z@\itemsep\z@\topsep\z@
         \setcounter{\@enumctr}{0}
         \def\makelabel##1{\hss\llap{\rm ##1}}
       }\fi}

\makeatother

\let\epsilon=\varepsilon
\def\({\big(}
\def\){\big)}

\def\0{\underline{0}}

\def\H{\mathscr H}

\def\B{\mathscr B_n}

\DeclareMathOperator{\Rad}{Rad}

\let\gdom\rhd
\let\gedom\unrhd

\def\m{\mathfrak m}

\def\m{\mathfrak m}

\def\s{\mathfrak s}

\def\t{\mathfrak t}
\def\u{\mathfrak u}
\def\v{\mathfrak v}

{\catcode`\|=\active
  \gdef\set#1{\mathinner{\lbrace\,{\mathcode`\|"8000%
                                   \let|\midvert #1}\,\rbrace}}
  \gdef\seT#1{\mathinner{\Big\lbrace\,{\mathcode`\|"8000%
                                   \let|\midverT #1}\,\Big\rbrace}}
}
\def\midvert{\egroup\mid\bgroup}
\def\midverT{\egroup\,\Big|\,\bgroup}

\def\Set[#1]#2|#3|{\Big\{\ #2\ \Big| \
           \vcenter{\hsize #1mm\centering #3}\Big\}}

\def\map#1#2{\,{:}\,#1\!\longrightarrow\!#2}

\title[  finite dimensional irreducible modules for  affine BMW
algebras]{On the classification of  finite dimensional irreducible modules for  affine BMW
algebras  }
\author{Hebing Rui}
\address{H.R. Department of Mathematics,  East China Normal
University, Shanghai, 200241, China} \email{hbrui@math.ecnu.edu.cn}

\thanks{The author was supported in part by NSFC and the Science and Technology Commission of
Shanghai Municipality 11XD1402200}

\begin{document}
\baselineskip15pt
\begin{abstract} In this paper, we classify the
finite dimensional irreducible modules for affine  BMW algebra over an algebraically closed field with arbitrary characteristic.
\end{abstract}
\sloppy \maketitle
%%%%%%%%%%%%%%%%%%%%%%%%%%%%%%%%%%%%%%%%%%%%%%%%%%%%%%%%%%%%%%%%%%%%%%%%%%
\section{Introduction }
%Throughout, we assume that $R$ is a commutative ring with multiplicative identity $1$ and units $\varrho, q$ and $q-q^{-1}$.

In \cite{HO:cycBMW}, Haering-Oldenburg introduced a class of associative algebras called affine Birman-Murakami-Wenzl (BMW for brevity)  algebras
in order to study knot invariants. These algebras, which can be considered as the affinization of BMW algebras in \cite{BW,  Mu}, had been studied extensively
by many authors  in
\cite{E, Enyang, GH, DRV, OR, G09,G12,RX, RS1, RS:bmw, RS2, RS4, W2, WY1, WY2} etc.

Recently,  Goodman~\cite{G12} studied  the cyclotomic
quotient of affine BMW algebras in $d$-semi-admissible case (see Definition~\ref{dsemi} for details).
This sets up the relationship between the representations of cyclotomic BMW algebras in general case and those for the cyclotomic BMW algebras in $\mathbf u$-admissible
case in  \cite{RX, RS1}.
Using the  results on the classification of irreducible modules of cyclotomic BMW algebras  in \cite{RX, RS1, Xi}, we get all finite dimensional irreducible modules for affine BMW algebras over an algebraically field $\kappa$ with arbitrary characteristic.

 In order to classify the finite  dimensional irreducible modules for affine BMW algebras over $\kappa$, we have to determine whether two irreducible modules for different cyclotomic BMW algebras are isomorphic as the modules for  the affine BMW algebra. For this, we need the result on the classification of finite dimensional irreducible modules for extended affine Hecke algebra $\hat\H_n$ of type $A_{n-1}$ as follows.

The first result on the classification of irreducible $\hat\H_n$-modules is due to Bernstein and Zelevinsky~\cite{BZ, Z}, who classified the irreducible $\hat\H_n$-modules over $\mathbb C$ when the defining parameter $q$ is not a root of unity. In this case, they used multisegments of length $n$  to index the complete set of non-isomorphic irreducible $\hat \H_n$-modules. In \cite{Ro}, Rogawski gave a different method to reprove Bernstien and Zelevinsky's result. Note that Kazhdan-Lusztig~\cite{KL} and Xi~\cite{Xinh} classified  the finite dimensional irreducible modules for affine Hecke algebras in any type. In particular, their results contain the case for extended affine Hecke algebras of type $A_{n-1}$.

On the other hand, any irreducible $\hat \H_n$-module over $\kappa$ can be realized as an irreducible module for an Ariki-Koike algebra~\cite{AK}. In the later case, its irreducible modules are indexed by Kleshchev multipartitions~\cite{Ari-osa}. In \cite{Va}, Vazirani gave the explicit relationship between the set of Kleshchev multi-partitions and the set of multi-segments when $q$ is not a root of unity.  If $q$ is a root of unity, the irreducible $\hat\H_n$-modules have been classified via aperiodic multisegments in \cite{Gin} (resp. \cite{AM}) over $\mathbb C$ (resp. over $\kappa$).
Further,  Ariki-Jacon-Lecouvey set up the explicit relationship between the set of Kleshchev multipartitions and the set of aperiodic multi-segments  in \cite[Theorem~6.2]{AJL} over $\kappa$. This is the result that we need when we classify the finite dimensional irreducible modules for
affine BMW algebras over $\kappa$.

Throughout, let   $\kappa$ be an algebraically closed field which contains non-zero elements
$q$,  $\varrho$, $\delta$ and a family of elements $\Omega=\{\omega_i\mid i\in \mathbb Z\}$
such that  $\delta=q-q^{-1}$ and $ \omega_0=1-\delta^{-1}(\varrho-\varrho^{-1})$.
Let $n$ be a positive integer  with $n\ge 2$.

\begin{Defn}\cite{HO:cycBMW}\label{affinebmw}  The affine  BMW algebra $\hat{\B}$
is the  unital associative $\kappa $--algebra generated  by $g_i,e_i,
x_1^{\pm 1}$, $1\le i\le n-1$ subject to the
following relations:
\begin{enumerate}
\item [(1)] $x_1 x_{1}^{-1}=x_{1}^{-1}x_1=1$ and $g_ig_i^{-1}=g^{-1}g_i=1$, for $1\le i\le n$,
\item  [(2)]  $g_ig_{i+1}g_i=g_{i+1}g_ig_{i+1}$, for $1\le i<n-1$,
\item  [(3)] $g_ig_j=g_jg_i$ if  $|i-j|>1$,
\item  [(4)] $x_1g_1x_1g_1=g_1x_1g_1x_1$, and $x_1 g_j=g_j x_1$ for $j\ge 2$,
\item  [(5)]  $e_i^2=\omega_0e_i$, for $1\le i<n$,
 \item  [(6)]  $e_1x_1^ae_1=\omega_a e_1$, for $a\in \mathbb Z^{>0}$,
    \item  [(7)]  $x_1g_j=g_jx_1$, for $2\le j\le n-1$,
    \item  [(8)]  $g_ie_j=e_jg_i$, and $e_ie_j=e_je_i$ if $|i-j|>1$,
\item  [(9)]  $e_ig_i=\varrho g_i=g_ie_i$, for $1\le i\le n-1$,
\item  [(10)]  $e_ig_{i\pm 1} e_i=\varrho e_i$, $e_{i}e_{i\pm 1} e_i=e_i$,
\item [(11)]  $g_ig_{i\pm 1} e_i=e_{i\pm 1} e_i$ and $e_i g_{i\pm 1} g_i=e_{i} e_{i\pm 1}$,
 \item  [(12)]  $g_i-g_{i}^{-1} =\delta (1-e_i)$, for $1\le i<n$,
\item   [(13)]  $e_1 x_1 g_1 x_1 g_1=e_1=g_1 y_1 g_1 x_1 e_1$.
\end{enumerate}
\end{Defn}

By Definition~\ref{affinebmw}, there is an anti-involution $\ast: \hat{\B}\rightarrow \hat{\B}$ which fixes $g_i, e_i$ and $x_1$, $1\le i\le n-1$.
Further, it is pointed in \cite[(2.1)]{GH} that Turaev~\cite{T}  has proved that $e_1 x_1^{-a} e_1=\omega_{-a} e_1 $ for $a\in \mathbb Z^{>0}$ and
$\omega_{-a}$ is a polynomial in $\omega_b$ for $b\in \mathbb Z^{>0}$. Therefore, $\omega_a$ is well-defined for all $a\in \mathbb Z$.

Goodman and Hauschild-Mosley~\cite{GH}constructed a basis for $\hat \B$ and showed that $\hat\B$ is of infinite dimension.
 In fact, Goodman and Hauschild-Mosley's  results~\cite{GH} are available over an integral domain.

 It is well-known that an affine Wenzl algebra in \cite{Na} can be considered  as a degenerate affine BMW algebra.
  Ariki, Mathas and Rui ~\cite{AMR} constructed an infinite dimensional irreducible modules for affine Wenzl algebra.
 Mimicking this  construction,  we know that $\hat\B$ has  infinite dimensional irreducible modules over a field.
 In other words, $\hat\B$ is not finitely generated over its center.  For the description of   the center of $\hat\B$, see \cite{DRV}.

The aim of this paper is to classify all finite dimensional  irreducible $\hat\B$-modules  over  $\kappa$.
Before we state our  main result, we need the notion of aperiodic multi-segments in \cite{AM}.

 Let $e$ be the smallest positive integer
such that $$1+q^2+q^4+ \cdots + q^{2(e-1)}=0$$ in $\kappa$\footnote{The current $q^2$ is $q$ in \cite{AM}.}. If there is no such  positive integer, then we set $e=\infty$. In other words, $e$ is the order of $q^2\in \kappa$.
Recall that a segment $\Delta$ of length $j=|\Delta|$ is a sequence of consecutive residues $[i, i+1, \ldots, i+j-1]$ where $i, i+1, \cdots, i+j-1\in \mathbb Z_{e}$. An multi-segment $\Delta$ is an unordered collection of segments $\Delta_i$ with length $\sum_i |\Delta_i|$.
Following \cite{AM}, we says that $\Delta$ is aperiodic if for every $j$, there is an $i\in \mathbb Z_e$ such that  $[i, i+1, \ldots, i+j-1]$
does not appear in $\Delta$.  Let  $\mathcal M_e^{n}$ be the   set of all aperiodic multi-segments with length $n$.
The following is the main result of this paper, which gives the classification of finite dimensional  irreducible $\hat\B$-modules over $\kappa$.

\begin{Theorem}\label{main1} Let $\hat\B$ be the affine BMW algebra over $\mathbb \kappa$.
\begin{enumerate} \item  Any finite dimensional irreducible  $\hat\B$-modules is of form $D^{f, \lambda}$ where $D^{f, \lambda}$, defined via the cellular basis of  some cyclotomic quotient $\mathscr B_{r, n}(\mathbf u)$ of $\hat \B$ in Theorem~\ref{Wcellular},   is an irreducible  $\mathscr B_{r, n}(\mathbf u)$-module   such that  \begin{enumerate} \item $0\le f\le \lfloor n/2\rfloor$ and $\lambda$ is a Kleshchev multipartition of $n-2f$ in the sense of \cite{AM}. Further, if $\omega_a=0$ for all $a\in \mathbb Z$ and if $2\mid n$, then $f\neq n/2$.
\item $\mathbf u$-admissible condition holds for $\mathscr B_{r, n}(\mathbf u)$ if $f>0$.\end{enumerate}
\item Let $D^{f, \lambda}$ (resp. $ D^{\ell, \mu}$) be the irreducible $\mathscr B_{r, n}(\mathbf u)$ (resp. $\mathscr B_{s, n}(\mathbf v)$)-module. Then   $D^{f, \lambda}\cong D^{\ell, \mu}$ as $\hat \B$-modules
 if and only if $f=\ell$ and the images of  $\lambda$ and $\mu$ under the map of \cite[Theorem~6.2]{AJL} are the same aperiodic multisegment in $\mathcal M_e^{n-2f}$ .
  \end{enumerate}
  \end{Theorem}

We remark that each aperiodic multisegment of length $n$ indexes an irreducible $\hat\B$-module on which $e_1$ acts trivially. This follows from Ariki-Mathas's result on the classification of irreducible $\hat\H_n$-modules  in \cite{AM}. However, we can not say that any pair $(f, \Delta)$ with $0<f<\lfloor n/2\rfloor$   and $\Delta\in \mathcal M_e^{n-2f}$ indexes an irreducible $\hat\B$-module. The reason is that each $\Delta\in \mathcal M_e^{n-2f}$ corresponds at least a Kleshchev multi-partition with respect to a family of parameters $u_1, u_2, \cdots, u_r\in \kappa^*$. However, we do not know whether the $\mathbf u$-admissible condition holds for $\mathscr B_{r, n}(\mathbf u)$.

 The content of this paper is organized as follows. We recall some of results on the representations of $\mathscr B_{r, n}(\mathbf u)  $ in section~2 and  prove Theorem~\ref{main1} in  section~3.

\section{Cyclotomic BMW algebras}
In this section, we recall some  results on the cyclotomic BMW algebra
over  $\kappa$ although some of them   hold over an integral domain.
Throughout, we assume  $r\in \mathbb Z$ with $r\ge 1$.

\begin{Defn}\label{cbmw}\cite{HO:cycBMW} Let $I$ be the two-sided ideal of $\hat \B$ generated by the cyclotomic polynomial
\begin{equation} \label{cycpoly} f(x_1)=(x_1-u_1)(x_1-u_2)\cdots (x_1-u_r),\end{equation}
where $u_i\in \kappa^*$, $1\le i\le r$.  The cyclotomic BMW  algebra $\mathscr B_{r, n}(\mathbf u)$ is the quotient algebra $\hat\B/ I$\footnote{In \cite{HO:cycBMW}, Haering-Oldenburg defined $\mathscr B_{r, n}(\mathbf u)$ without assuming $u_i\in \kappa^*$, $1\le i\le r$.}.
\end{Defn}

\begin{Remark}  When $r=1$,  $\mathscr B_{r, n}(\mathbf u)$ is the usual BMW algebra, which was introduced by Birman-Wenzl~\cite{BW} and independently by Murakami~\cite{Mu}.\end{Remark}

It  is known that  $\mathscr B_{r, n}(\mathbf u)$  can be used to study the finite dimensional irreducible $\hat\B$-modules over $\kappa$.
Pick a finite dimensional irreducible $\hat\B$-module $M$ over $\kappa$. Let $f(x_1)$ be the characteristic polynomial of $x_1$ with respect to $M$. Then
$M$ has to be an irreducible $\mathscr B_{r, n}(\mathbf u)$-module where  $\mathscr B_{r, n}(\mathbf u)=\hat \B/ I$
and  $I$ is the two-sided ideal of $\hat\B$ generated by $f(x_1)$.
Since $\kappa$ is an algebraically close field,  $f(x_1)$ is given in (\ref{cycpoly}) for  some   $u_1, u_2, \cdots, u_r\in \kappa$. Further,  $u_i\in \kappa^*$
 since $x_1$ is invertible in $\mathscr B_{r, n}(\mathbf u)$. Therefore, we  will get  all finite dimensional irreducible $\hat\B$-modules over $\kappa$ if we classify the irreducible
 $\mathscr B_{r, n}(\mathbf u)$-modules for all $\mathbf u\in (\kappa^*)^r$ and $r\ge 1$.

\begin{Defn}\label{dsemi}\cite{G12} We say that the $d$-semi-admissible condition holds for  $\mathscr B_{r, n}(\mathbf u)$
if $d$ is the minimal integer such that  $\{e_1, e_1x_1, \cdots, e_1 x_1^{d}\}$ is linear dependent in $\mathscr B_{r, 2}(\mathbf u)$. \end{Defn}

Obviously,  $0\le d\le r$. We have  $e_1=0$ if $d=0$.  In this case, there is no restriction on $u_i$'s. Further,  $\mathscr B_{r, n}(\mathbf u)$ is the Ariki-Koike algebra $\H_{r, n}$~\cite{AK} whose simple modules have been classified in \cite{Ari-osa}.

If $d=r$,  then the $d$-semi-admissible condition is the $\mathbf u$-admissible condition in \cite{RX} or  admissible conditions in \cite{WY1}. In particular,  $\mathbf u$-admissible condition always holds  if $e_1\neq 0$ and $r=1$.

In $\mathbf u$-admissible case, we have~\cite{RX}

 \begin{equation} \varrho^{- 1}= \alpha \prod_{\ell=1}^r
u_{\ell}, \text{ and } \omega_{a} = \sum\limits_{j=1}^r
u_{j}^{a}\gamma_{j}, \forall a\in \mathbb Z,\end{equation}  where
\begin{itemize} \item [(1)] $\gamma_{i} = (\gamma_r(u_i) + \delta^{-1}\varrho(u_{i}^{2}
- 1) \prod \limits_{j\neq i}u_{j})\prod\limits_{j\neq
i}\frac{u_{i}u_{j} - 1}{u_{i} - u_{j}}$, and $\gamma_r(z)$ is $1$ (resp.  $-z$)  if $2\nmid r$ (resp.  otherwise).

\item [(2)] $\alpha\in \{1, -1\}$ if $2\nmid r$ and $\alpha\in \{q^{-1}, -q\}$, otherwise.
\item[(3)]
 $\omega_{0}=
\delta^{-1}\varrho(\prod\limits_{\ell=1}^r u_{\ell}^2 - 1) +
1-\frac{(-1)^r +1}{2}\alpha^{-1}\varrho^{-1} $.
\end{itemize}

We have the following result, which will be used when we prove Theorem~\ref{main1}.

 \begin{Lemma}\label{zero} Suppose  $\mathbf u$-admissible condition holds for $\mathscr B_{r, 2}(\mathbf u)$. We have $\omega_i\neq 0$ for some $i, 0\le i\le r-1$ if
 there is a $j\in \mathbb Z$ such that $\omega_j\neq 0$.
 \end{Lemma}
\begin{proof}
This follows from  Definitions~\ref{cbmw} and~\ref{affinebmw}(6).
\end{proof}

If the $\mathbf u$-admissible condition holds, then
 $\mathscr B_{r, n}(\mathbf u)$ is (weakly) cellular in the sense of
\cite{GL} as follows.

\begin{Defn}  \cite{GL}\label{GL} Assume that $R$ is a commutative ring with the multiplicative identity $1$.
    Let  $A$ be  an $R$--algebra.
    Fix a partially ordered set $\Lambda=(\Lambda,\gedom)$ and for each
    $\lambda\in\Lambda$ let $T(\lambda)$ be a finite set. Finally,
    fix $\m_{\s\t}\in A$ for all
    $\lambda\in\Lambda$ and $\s,\t\in T(\lambda)$.
    Then the triple $(\Lambda,T,C)$ is a \textsf{cell datum} for $A$ if:
    \begin{enumerate}
    \item $\mathcal M=\set{\m_{\s\t}|\lambda\in\Lambda\text{ and }\s,\t\in
        T(\lambda)}$ is an $R$--basis for $A$;
    \item the $R$--linear map $*\map AA$ determined by
        $(\m_{\s\t})^*=\m_{\t\s}$, for all
        $\lambda\in\Lambda$ and all $\s,\t\in T(\lambda)$ is an
        anti--isomorphism of $A$;
    \item for all $\lambda\in\Lambda$, $\s\in T(\lambda)$ and $a\in A$
        there exist scalars $r_{\t\u}(a)\in R$ such that
        $$\m_{\s\t} a
            =\sum_{\u\in T(\lambda)}r_{\t\u}(a)\m_{\s\u}
                     \pmod{A^{\gdom\lambda}},$$
            where
    $A^{\gdom\lambda}=R\text{--span}%
      \set{\m_{\u\v}|\mu\gdom\lambda\text{ and }\u,\v\in T(\mu)}$.
     Furthermore, each scalar $r_{\t\u}(a)$ is independent of $\s$. \end{enumerate}
     An algebra $A$ is a \textsf{cellular algebra} if it has
    a cell datum and in this case we call $\mathcal M$
    a cellular basis of $A$.
\end{Defn}

The notion of weakly cellular algebras in \cite{G} is
obtained from Definition~\ref{GL} by using $$(\m_{\s\t})^*\equiv
\m_{\t\s} \pmod {A^{\gdom\lambda}}$$ instead of
$(\m_{\s\t})^*=\m_{\t\s}$. Note that both cellular algebras and
weakly cellular algebras are standardly based algebras in the sense
of \cite{DR}. From this, one can see that  cellular algebras and
weakly cellular algebras share the similar results on representation
theory. For this reason, both cellular algebras  and weakly
cellular algebras will be called cellular algebras later on.

Now, we briefly recall the representation theory of cellular
algebras over a field in ~\cite{GL}. We remark that all modules
considered in this paper are right modules.

Every irreducible $A$--module arises in a unique way as the simple
head of some cell module.
For each $\lambda\in\Lambda$ fix $\s \in
T(\lambda)$ and let $$\m_{\t}
       =\m_{\s\t}+ A^{\rhd \lambda}.$$
       The cell module $S^{\lambda}$ of $A$ with respect to $\lambda$ can be considered as
        the free $R$--module with basis $\set{\m_{\t}|\t\in
T(\lambda)}$. The cell module $S^\lambda$ comes equipped with
a natural symmetric bilinear form $\phi_{\lambda}$ which is
determined by the equation
$$\m_{\s\t}
           \m_{\t'\s}
 \equiv\phi_{\lambda}\big(\m_{\t},
               \m_{\t'}\big)\cdot
        \m_{\s\s}\pmod{ A^{\rhd \lambda}}.$$
The bilinear form $\phi_{\lambda}$ is $A$--invariant in the sense
that $$\phi_{\lambda}(xa,y)=\phi_{\lambda}(x,ya^*), \text{for
$x,y\in S^\lambda$ and $a\in A$.}$$ Consequently,
$$\Rad S^\lambda
   =\set{x\in S^\lambda|\phi_{\lambda}(x,y)=0\text{ for all }
                          y\in S^\lambda}$$
is an $A$--submodule of $S^\lambda$ and
$D^\lambda=S^\lambda/\Rad S^\lambda$ is either zero or
absolutely irreducible.

Graham and Lehrer~\cite{GL} have proved that all non-zero
 $D^\lambda$  consist of a complete set
of pairwise non-isomorphic irreducible $A$-modules.  This gives a useful method to classify the
irreducible modules for cellular algebras.

Recall that a \textsf{composition} $\lambda=(\lambda_1,\lambda_2,\dots)$ of $m$ is a sequence of non-negative integers
  with
$|\lambda|=\sum_i\lambda_i=m$.  If $\lambda$ is weakly decreasing, then $\lambda$ is called a partition. Similarly, an
{$r$-partition} of $m$ is an ordered $r$-tuple
$\lambda=(\lambda^{(1)},\dots,\lambda^{(r)})$ of partitions
$\lambda^{(s)}$ with $1\le s\le r$ such that
$|\lambda|=\sum_{i=1}^r|\lambda^{(i)}|=m$. Let
$\Lambda_r^+(n)$ be the set of all $r$-partitions of $n$.
We say
that $\mu$ dominances $\lambda$ and write
$\lambda\trianglelefteq\mu$ if  $$\sum_{j=1}^{i-1} |\lambda^{(j)}|
+\sum_{k=1}^l \lambda_k^{(i)} \le \sum_{j=1}^{i-1}
|\mu^{(j)}|+\sum_{k=1}^l \mu_k^{(i)} $$ for $1\le i\le r$ and
$l\ge 0$. So, $(\Lambda_r^+(n), \unlhd)$ is a poset. If
$\lambda\trianglelefteq\mu$ and $\lambda\neq \mu$, we write
$\lambda\vartriangleleft\mu$.
 Let \begin{equation} \label{posetc} \Lambda_{r,
n}=\{(k, \lambda)\mid 0\le k\le \lfloor n/2\rfloor, \lambda\in
\Lambda_r^+(n-2k)\}.\end{equation}  Then  $\Lambda_{r, n}$ is a poset with
$\unrhd$ as the  partial order on it.
 More explicitly, $(k, \lambda)\unrhd(\ell, \mu) $
for $(k, \lambda), (\ell, \mu)\in \Lambda_{r, n}$ if either
$k>\ell$ in the usual sense or $k=\ell$ and $\lambda\unrhd\mu$. Here
$\unrhd$ is the dominance order defined on  $\Lambda_r^+(n-2k)$.

The following theorem is well-known. See   \cite{G,  WY2} for another description of cellular basis for $\mathscr B_{r, n}(\mathbf u)$.

\begin{Theorem}\cite{RX} \label{Wcellular} Suppose that the  $\mathbf u$-admissible condition holds for $\mathscr B_{r, n}(\mathbf u)$.
 Then
$$\mathscr C=\bigcup_{(f, \lambda)\in \Lambda_{r,
n}}\set{C_{\s\t}|
            \s ,\t \in T(f,\lambda)
              }$$
is a weakly cellular basis of  $\mathscr B_{r, n}(\mathbf u)$ over the poset $\Lambda_{r, n}$. In this case, the required
 $\kappa$-linear  anti-involution
 on $\mathscr B_{r, n}(\mathbf u)$ is $\ast$, which  fixes $g_i, e_i$ and $x_1, 1\le i\le n-1$. In particular, the rank of $\mathscr B_{r, n}(\mathbf u)$ is $r^n (2n-1)!!$.
\end{Theorem}
In this paper, we do not need
 the explicit definition of $C_{\s\t}$ in  \cite[4.17]{RX}.     What we will need is  some properties of cell modules $S^{f, \lambda}$, $(f, \lambda)\in \Lambda_{r, n}$
  for $\mathscr B_{r, n}(\mathbf u)$ with respect to the cellular basis in Theorem~\ref{Wcellular}.   Let $\phi_{f, \lambda}$ be the invariant form on the cell
module $S^{f, \lambda}$ with respect to $\lambda\in \Lambda_{r, n}$.

We have  $\mathscr B_{r, n}(\mathbf u)/I\cong \H_{r, n}(\mathbf u)$, where $\H_{r, n}(\mathbf u)$
is the Ariki-Koike algebra~\cite{AK} and $I$ is the two-sided deal of $\mathscr B_{r, n}(\mathbf u)$ generated by the cyclotomic polynomial
$f(x_1)$ in (\ref{cycpoly}). The image of the cellular basis of $\mathscr B_{r, n}(\mathbf u)$ in Theorem~\ref{Wcellular} is the cellular basis of $\H_{r, n}(\mathbf u)$ in \cite{DJM}.
The corresponding cell module of $\H_{r, n}(\mathbf u)$ with respect to $\lambda\in \Lambda_r^+(n)$ is denoted by $S^\lambda$. Let $\phi_\lambda$ be the invariant form on $S^\lambda$.

%The following result has been proved in \cite{RX}.

\begin{Prop}~\label{cbmw-equiv}\cite[5.2]{RX}  Suppose that the $\mathbf u$-admissible condition holds for   $\mathscr B_{r, n}(\mathbf u)$ over  $\kappa$\footnote{In \cite{RX}, $\kappa$ is an arbitrary field. }.
 Assume that $(f, \lambda)\in \Lambda_{r, n}$.
\begin{enumerate} \item
 If $f\neq n/2$, then  $\phi_{f, \lambda}\neq 0$ if and only if $\phi_\lambda\neq 0$.
\item  If $\omega_a\neq 0$ for  some non-negative integer $a\le r-1$, then $\phi_{n/2, 0}\neq 0$.
\item  If  $\omega_a=0$ for   all  non-negative integers $a\le r-1$, then  $\phi_{n/2, 0}= 0$. \end{enumerate}
\end{Prop}

Note that $\phi_\lambda\neq 0$ if and only if $D^\lambda\neq 0$ for $\H_{r, n-2f}(\mathbf u)$.
By \cite{Ari-osa},
$\phi_{\lambda}\neq 0$ if and only if $\lambda$ is a \textsf{ Kleshchev multipartition } in   the sense of the Definition in \cite[p605]{AM}.
So, the irreducible $\mathscr B_{r, n}(\mathbf u)$-modules are classified via Proposition~\ref{cbmw-equiv}. More explicitly, we have the following result which can be found in \cite{Xi} for $r=1$ and
~\cite{RX} for  $r\ge 2$. We remark that the $\mathbf u$ admissible condition always holds for $r=1$ and $e_1\neq 0$.

\begin{Theorem}~\label{cbmw-simple}\cite{RX, Xi}  Suppose that $\mathbf u$-admissible condition holds for  $\mathscr B_{r, n}(\mathbf u)$  over $\kappa$.

\begin{enumerate} \item If either $\omega_a\neq 0$ for  some non-negative integer $a\le r-1$ or  $\omega_a=0$ for   all  non-negative integers $a\le r-1$ and $2\nmid n$, then the
irreducible  $\mathscr B_{r, n}(\mathbf u)$-modules are indexed by $(f, \lambda)$ with $0\le f\le \lfloor n/2\rfloor$
and $\lambda$'s are Klechshev multipartions of $n-2f$.
\item  If   $\omega_a=0$ for   all  non-negative integers $a\le r-1$ and $2\mid n$, then the
irreducible  $\mathscr B_{r, n}(\mathbf u)$-modules are indexed by $(f, \lambda)$ with $0\le f< \lfloor n/2\rfloor$
and $\lambda$'s are   Klechshev multipartions of $n-2f$.
 \end{enumerate}
\end{Theorem}

At the end of this section, we recall Goodman's result for  $0<d<r$ in \cite{G12}. In this case, $d$ is the  minimal integer such that  $\{e_1, e_1x_1, \cdots, e_1 x_1^{d}\}$  is linear dependent  in $\mathscr B_{r, 2}(\mathbf u)$. Goodman~\cite{G12} showed that there is a polynomial $g(x_1)\in \kappa [x_1]$ with $\text{deg.} g(x_1)=d$ such that $e_1g(x_1)=0$ and $e_1 h(x_1)\neq 0$ in $\mathscr B_{r, 2}(\mathbf u)$ for any polynomial $h(x_1)\in \kappa[x_1]$ with  $\text{deg.} h(x_1)<d$.
Further, since $e_1f(x_1)=0$ in $\mathscr B_{r, 2}(\mathbf u)$, it is not difficult to see that $g(x_1)\mid f(x_1)$.
So, write $$g(x_1)=(x_1-v_1)(x_1-v_2)\cdots (x_1-v_d),$$
where $\{v_1, v_2, \cdots, v_d\}\subset \{
u_1, u_2, \cdots, u_r\}$ such that $\mathbf v$-admissible condition holds in $\mathscr B_{d, n}(\mathbf v)$. Let $\langle e_1\rangle_r$ (resp. $\langle e_1\rangle_d$) be the two-sided ideal of
$\mathscr B_{r, n}(\mathbf u)$ (resp. $\mathscr B_{d, n}(\mathbf v)$) generated by $e_1$.

\begin{Theorem}\cite[5.11]{G12}  There is an algebraically epimorphism $\theta:\mathscr B_{r, n}(\mathbf u)\twoheadrightarrow\mathscr B_{d, n}(\mathbf v)$ such that the restriction of
$\theta$ on $\langle e_1\rangle_r$ gives rise to  an isomorphism between $\langle e_1\rangle_r$ and   $\langle e_1\rangle_d$.
\end{Theorem}

Since  $\mathbf v$-admissible conditions hold in $\mathscr B_{d, n}(\mathbf v)$, $\langle e_1\rangle_d$ is  cellular with a  basis which is
given in  Theorem~\ref{Wcellular} for $\mathscr B_{d, n}(\mathbf v)$ with respect to the poset which consists of all pairs $(f, \lambda)\in \Lambda_{d, n}$ such that $f\ge1$.
Via the isomorphism $\theta$, Goodman~\cite{GL} lifted the cellular basis of $\langle e_1\rangle_d$ to get the corresponding cellular basis of  $ \langle e_1\rangle_r$.
Using  the epimorphism
$\pi: \mathscr B_{r, n}(\mathbf u)\twoheadrightarrow \H_{r, n}(\mathbf u)$,   Goodman~\cite{G12} showed the following result.

\begin{Theorem}\cite[Theorem~6.4]{G12} \label{dsemi} Suppose that the $d$-semi-admissible condition holds for $\mathscr B_{r, n}(\mathbf u)$. Then  $\mathscr B_{r, n}(\mathbf u)$ is (weakly) cellular over the poset
$$\tilde \Lambda_{r, n}=\cup_{1\le f\le \lfloor n/2\rfloor} \{(f, \lambda)\mid (f,\lambda)\in \Lambda_d^+(n-2f)\}\cup\{(0, \lambda)\mid \lambda\in \Lambda_r^+(n)\} $$ in the sense of Definition~\ref{GL}. Further, $\dim_\kappa \mathscr B_{r, n}(\mathbf u)=d^n (2n-1)!! +r^n n! -d^n n!$.\end{Theorem}

We remark that  $(f, \lambda)\le (\ell, \mu)$ for  $(f, \lambda), (\ell, \mu)\in \tilde  \Lambda_{r, n}$ if either  $f<\ell $ or $f=\ell$ and $\lambda\unlhd \mu$ where $\unlhd$ is the dominance order on $\Lambda_d^+(n-2f)$ (resp. $\Lambda_r^+(n)$) provided $f>0$ (resp. $f=0$).

For each $(f, \lambda)\in \Lambda_{d, n}$ with $f\ge 1$, let  $S^{f, \lambda}$ (resp.  $D^{f, \lambda}$) be the cell (resp. irreducible ) module of   $\mathscr B_{d, n}(\mathbf v)$ with respect the cellular basis in Theorem~\ref{Wcellular}. Then  $S^{f, \lambda}$ (resp.  $D^{f, \lambda}$) can be  considered as the corresponding cell  (resp. irreducible ) module  of
$\mathscr B_{r, n}(\mathbf u)$ with respect to  $(f, \lambda)\in \tilde \Lambda_{r, n}$ such that $f>0$. Therefore, we can always assume that $\mathbf u$-admissible conditions holds when we discuss the irreducible module  $D^{f, \lambda}$ for $f>0$. This is the reason why we add $\mathbf u$-admissible  condition in Theorem~\ref{main1}(a)(ii).

\section{Proof of Theorem~\ref{main1}  }
In this section, we prove Theorem~\ref{main1}, which gives the classification of finite dimensional irreducible $\hat \B$-modules over $\kappa$.

\begin{Lemma} \label{ebe} Suppose $n> 2$. If $\omega_0\neq 0$, we define $e=\omega_0^{-1} e_{n-1}$. Otherwise,
we define $e=\rho^{-1} e_{n-1} g_{n-2}$. Then $e^2=e$ and  $e \hat \B e= \hat{\mathscr B}_{n-2}e\cong \hat{\mathscr B}_{n-2}$   as $\kappa$-algebras.
\end{Lemma}

\begin{proof} It follows from  Definition~\ref{affinebmw}(5)(10) that   $e^2=e$. By \cite[3.17,~3.20]{GH},   $e_{n-1} \hat {\mathscr B}_{n-1} e_{n-1}=e_{n-1} \hat{\mathscr B}_{n-2}$
and $ \hat {\mathscr B}_{n} e_{n-1} =\hat {\mathscr B}_{n-1} e_{n-1}$.  Therefore, $e_{n-1} \hat \B e_{n-1}=e_{n-1} \hat{\mathscr B}_{n-2}$.
Now, everything follows since  $g_{n-2}$ is invertible. We remark that the required  isomorphism from $\hat{\mathscr B}_{n-2}$ to  $\hat{\mathscr B}_{n-2}e$ sending $x$ to $xe$ for all
$x\in  \hat{\mathscr B}_{n-2}$.  One can verify the injectivity of this homomorphism by using   the result on the basis  of $\hat \B$ in \cite{GH}.
\end{proof}

Let $\hat\B$-mod be the category of finite dimensional right $\hat \B$-modules over $\kappa$. By Lemma~\ref{ebe},  we have
the functor $\mathfrak F:  \text{$\hat\B$-mod}\rightarrow  \text {$ \hat{\mathscr B}_{n-2}$-mod}$ such that \begin{equation} \label{func} \mathfrak F(M)=M e\end{equation}
for any object $M\in  \hat\B\text{-mod}$. Further, if $M$ is a $\mathscr B_{r, n}(\mathbf u)$-module  and if there is an epimorphism  $\phi: \hat \B\twoheadrightarrow \mathscr B_{r, n}(\mathbf u)$,   then $\mathfrak F(M)$ is the same as $\mathcal F(M)$ where
$\mathcal F$ is the exact functor from $\mathscr B_{r, n}(\mathbf u)$-mod to $\mathscr B_{r, n-2}(\mathbf u)$-module.     However, by Theorem~\ref{dsemi} and the statements below Theorem~\ref{dsemi}, we can always assume the $\mathbf u$-admissible condition holds when we discuss $S^{f, \lambda}$ and $D^{f, \lambda}$ for $f>0$. In this case, by ~\cite[Sect. 5]{RS1}, we  have  \begin{equation} \mathcal F(S^{f, \lambda})=S^{f-1, \lambda}\text{ and
 $\mathcal F(D^{f, \lambda})=D^{f-1, \lambda}$ }.\end{equation}  Note that $D^{f, \lambda}\neq 0$ if and only of $D^{ \lambda}\neq 0$ (see   Proposition~\ref{cbmw-equiv}).
Further,
$$\mathcal F(S^{0, \lambda})=\mathcal F(D^{0, \lambda})=0$$
no matter whether the $\mathbf u$-admissible condition holds for  $\mathscr B_{r, n}(\mathbf u)$.

\begin{Lemma} \label{key1} Suppose $(f, \lambda)\in \Lambda_{r, n}$ (resp.  $(\ell, \mu)\in \Lambda_{s, n}$) such that
 $D^{f, \lambda}\neq 0$ (resp.  $D^{\ell, \mu}\neq 0$) as $\mathscr B_{r, n}(\mathbf u)$-module  (resp. $\mathscr B_{s, n}(\mathbf v)$-module).
 If both $\mathscr B_{r, n}(\mathbf u)$ and $\mathscr B_{s, n}(\mathbf v)$ are images of $\hat\B$ such that
$D^{f, \lambda}\cong D^{\ell, \mu}$ as $\hat\B$-modules, then $f=\ell$ and $D^{\lambda}\cong D^{ \mu}$ as $\hat \H_{n-2f}$-modules.
\end{Lemma}
 \begin{proof} If $f\neq \ell$, we can assume that $f\ge \ell+1$ without loss of any generality.
 By Theorem~\ref{dsemi}, we can always assume that $\mathbf u$-admissible (resp. $\mathbf v$-admissible ) condition holds (resp.  if $\ell\neq 0$).

 Applying the functor $\mathfrak F$ on both
 $D^{f, \lambda}$ and $ D^{\ell, \mu}$ repeatedly yields $D^{f-\ell-1, \lambda}= 0$ as $\hat{\mathscr {B}}_{n-2\ell-2}$-modules. By Proposition~\ref{cbmw-equiv},
 $D^{f-\ell-1, \lambda}\neq 0$, a contradiction. So, $f=\ell$.  Applying the functor $\mathfrak F$ on both
 $D^{f, \lambda}$ and $ D^{f, \mu}$ yields $D^{0, \lambda}\cong D^{0, \mu}$ as $\hat{\mathscr B}_{n-2f}$-modules.
 In other words,  $D^{0, \lambda}\cong D^{0, \mu}$ as $\hat\H_{n-2f}$-modules.
Note that $D^{0, \lambda}$ (resp. $D^{0, \mu}$)  can be identified with  $D^\lambda$ (resp. $D^\mu$) as $ \H_{r, n-2f}(\mathbf u)$
(resp. $ \H_{s, n-2f}(\mathbf v)$)-module. Now, everything follows.
 \end{proof}

\begin{Lemma}\label{key3} Suppose   $(f, \lambda)\in \Lambda_{r, n}$ (resp.  $(f, \mu)\in \Lambda_{s, n}$) such that
 $D^{f, \lambda}\neq 0$ (resp.  $D^{f, \mu}\neq 0$) as $\mathscr B_{r, n}(\mathbf u)$-module  (resp. $\mathscr B_{s, n}(\mathbf v)$-module).
 If both $\mathscr B_{r, n}(\mathbf u)$ and $\mathscr B_{s, n}(\mathbf v)$ are images of $\hat\B$ and if  $D^{ \lambda}\cong D^{ \mu}$ as $\hat \H_{n-2f}$-modules, then
$D^{f, \lambda}\cong D^{f, \mu}$ as $\hat\B$-modules.
\end{Lemma}

\begin{proof}  First, we can assume $f\neq 0$. Otherwise, there is nothing to be proved.
 By assumption,  $D^{f, \lambda}$ (resp. $D^{f, \mu}$) is the irreducible $\mathscr B_{r, n}(\mathbf u)$-module (resp.  $\mathscr B_{s, n}(\mathbf v)$-module)
with respect to  $(f, \lambda)\in \Lambda_{r, n}$ (resp.  $(f, \mu) \in \Lambda_{s,n}$).
Suppose
$$\mathscr B_{r, n}=\hat{\mathscr B}_n/I \text{ and } \mathscr B_{s, n}=\hat{\mathscr B}_n/J,$$ where $I$ (resp. $J$)  is the two-sided ideal of $\hat{\mathscr B}_n$ generated by
$f(x_1)$ (resp. $g(x_1)$ ) and
$$ \begin{aligned} f(x_1) & =(x_1-u_1)(x_1-u_2)\cdots (x_1-u_r), \\  g(x_1)&=(x_1-v_1)(x_1-u_2)\cdots (x_1-v_s).\\
\end{aligned}$$
Let $h(x)=[f(x_1), f(x_2)]$ be the least common multiple of $f(x_1)$ and $g(x_1)$. Let $\mathscr B_{t, n}=\hat{\mathscr B}_n/K$ where $K$
  is the two-sided ideal of $\hat{\mathscr B}_n$ generated by $h(x_1)$. Then there are two algebraical epimorphisms:
  $$ \phi: \mathscr B_{t, n}\twoheadrightarrow  \mathscr B_{r, n}(\mathbf u), \text{ and }  \psi: \mathscr B_{t, n}\twoheadrightarrow \mathscr B_{s, n}(\mathbf v)$$
such that  $\phi$ (resp.  $\psi$) sends generators $e_i, g_i, x_1\in \mathscr B_{t, n} $ to the corresponding generators $e_i, g_i, x_1$ in $\mathscr B_{r, n}(\mathbf u)$ (resp. $\mathscr B_{s, n}(\mathbf v)$), $1\le i\le n-1$.

In particular, by Lemma~\ref{key1} and Theorem~\ref{dsemi},  the irreducible $\mathscr B_{r, n}(\mathbf u) $-module (resp.  $\mathscr B_{s, n}(\mathbf v) $-module)   $D^{f, \lambda}$ (resp. $D^{f, \mu}$)
    has to be  the irreducible  $\mathscr B_{t, n}$-module
$D^{f, \alpha}$ (resp. $D^{f, \beta}$)  for some multi-partition $\alpha$ (resp. $\beta$) such that $D^\alpha\cong D^\lambda$ and  $D^\beta\cong D^\mu$ as $\hat\H_{n-2f}$-modules.
Further, by the arguments below Theorem~\ref{dsemi} and the results on the representation theory for cyclotomic BMW algebras $\mathscr B_{t, n}$,
we have that both $D^\alpha$ and $D^\beta$ are irreducible modules for the same Ariki-Koike algebra.
Since we are assuming that $ D^\lambda\cong D^\mu$ as $\hat\H_{n-2f}$-modules, we have $D^\alpha\cong D^{\beta}$, forcing  $\alpha=\beta$. So, $D^{f, \lambda}\cong D^{f, \alpha}\cong  D^{f, \mu}$ as $\hat\B$-modules, proving the result.
\end{proof}

\textsf{Proof of Theorem~\ref{main1}:} Let $D^\lambda$ be the irreducible modules for $\H_{r, n}(\mathbf u)$.
In \cite{Ari-osa}, Ariki has proved that $D^\lambda\neq 0$ if and only if $\lambda$ is Kleshchev in the sense of Definition in \cite[p605]{AM}. On the other hand, simple modules for affine Hecke algebra $\hat \H_n$ can be labeled by aperiodic multisegments~\cite{AM}. In \cite{AJL}, Ariki, Jacon and Lecouvey set up the explicit relationship
between the set of Kleshchev multipartitions and the set of aperiodic multisegments in \cite[Theorem~6.2]{AJL}.  In other words, if  $D^\lambda$ and $D^\mu$ are irreducible modules for different
Ariki-Koike algebras with respect to the Kleshchev multi-partitions $\lambda$ and $\mu$, then $D^\lambda\cong D^\mu$ as irreducible modules for extended affine Hecke algebra  if and only if the images of $\lambda$ and $\mu$ with respect to (different) map in  \cite[Theorem~6.2]{AJL} are the same aperiodic
multi-segment with length $n-2f$. Fuether, when $f>0$, we have to assume the $\mathbf u$-admissible condition hold. However, when $f=0$, we do not need this assumption.
 Now, everything follows from Lemmas~\ref{key1}-\ref{key3} and ~\ref{zero}.\qed

We close the paper by giving  the following remark.
 
 \begin{Remark} We can classify the finite dimensional irreducible modules for affine Wenzl algebra over an algebraically closed field $\kappa$. In this case, we have to use the results for degenerate affine Hecke algebra of type $A_{n-1}$ instead of those for $\hat\H_n$.  We leave the details to the reader.
\end{Remark}

\providecommand{\bysame}{\leavevmode ---\ } \providecommand{\og}{``}
\providecommand{\fg}{''} \providecommand{\smfandname}{and}
\providecommand{\smfedsname}{\'eds.}
\providecommand{\smfedname}{\'ed.}
\providecommand{\smfmastersthesisname}{M\'emoire}
\providecommand{\smfphdthesisname}{Th\`ese}

\end{document}